\documentclass[reqno,11pt]{amsart}

\usepackage{amsmath,amsfonts,amssymb,graphicx,amsthm,enumerate,url}
\usepackage[noadjust]{cite}
\usepackage{stmaryrd}

\usepackage{xifthen,xcolor,tikz,setspace}
\usetikzlibrary{decorations.pathmorphing,patterns,shapes,calc}
\usepackage[noadjust]{cite}
\definecolor{refkey}{gray}{.75}
\definecolor{labelkey}{gray}{.5} 

\usepackage[colorlinks=true]{hyperref}
\colorlet{DarkGreen}{green!50!black}
\colorlet{DarkGray}{gray!60!black}

\numberwithin{equation}{section}

\renewcommand{\restriction}{\mathord{\upharpoonright}}
\renewcommand{\epsilon}{\varepsilon}
\newcommand{\given}{\;\big|\;}

\usepackage{color}
 \definecolor{refkey}{gray}{.5}
 \definecolor{labelkey}{gray}{.5}
\definecolor{light}{gray}{.9}

\newtheorem{maintheorem}{Theorem}

\newtheorem{theorem}{Theorem}[section]
\newtheorem*{theorem*}{Theorem}
\newtheorem{lemma}[theorem]{Lemma}

\newtheorem{proposition}[theorem]{Proposition}

\newtheorem{fact}[theorem]{Fact}
\newtheorem{corollary}[theorem]{Corollary}

\theoremstyle{definition}{

\newtheorem{definition}[theorem]{Definition}
\newtheorem*{definition*}{Definition}

\newtheorem{question}[theorem]{Question}
\newtheorem{remark}[theorem]{Remark}
}

\newcommand{\cA}{\ensuremath{\mathcal A}}
\newcommand{\cB}{\ensuremath{\mathcal B}}

\newcommand{\cF}{\ensuremath{\mathcal F}}
\newcommand{\cG}{\ensuremath{\mathcal G}}

\newcommand{\cJ}{\ensuremath{\mathcal J}}

\newcommand{\cM}{\ensuremath{\mathcal M}}

\newcommand{\cR}{\ensuremath{\mathcal R}}

\newcommand{\cY}{\ensuremath{\mathcal Y}}

\renewcommand{\epsilon}{\varepsilon}

\newcommand{\bP}{{\mathbf{P}}}

\makeatletter
\newcommand{\superimpose}[2]{%
  {\ooalign{$#1\@firstoftwo#2$\cr\hfil$#1\@secondoftwo#2$\hfil\cr}}}
\makeatother

\begin{document}

\title[Zero-temperature dynamics in the dilute Curie--Weiss model]{Zero-temperature dynamics in the dilute Curie--Weiss model}

\author{Reza Gheissari}
\address{R.\ Gheissari\hfill\break
Courant Institute\\ New York University\\
251 Mercer Street\\ New York, NY 10012, USA.}
\email{reza@cims.nyu.edu}

\author{Charles M.\ Newman}
\address{C.\ M.\ Newman\hfill\break
Courant Institute\\ New York University\\
251 Mercer Street\\ New York, NY 10012, USA.}
\email{newman@cims.nyu.edu}

\author{Daniel L.\ Stein}
\address{D.\ L.\ Stein\hfill\break
Courant Institute\\ New York University\\
251 Mercer Street\\ New York, NY 10012, USA.}
\email{daniel.stein@nyu.edu}

\begin{abstract}We consider the Ising model on a dense Erd\H{o}s--R\'enyi random graph, $\mathcal G(N,p)$, with $p>0$ fixed---equivalently, a disordered Curie--Weiss Ising model with $\mbox{Ber}(p)$ couplings---at zero temperature. The disorder may induce local energy minima in addition to the two uniform ground states. 
In this paper we prove that, starting from a typical initial configuration, the zero-temperature dynamics avoids all such local minima and absorbs into a predetermined one of the two uniform ground states. We relate this to the local MINCUT problem on dense random graphs; namely with high probability, the greedy search for a local MINCUT of $\mathcal G(N,p)$ with $p>0$ fixed, started from a uniform random partition, fails to find a non-trivial cut. In contrast, in the disordered Curie--Weiss model with heavy-tailed couplings, we demonstrate that zero-temperature dynamics has positive probability of absorbing in a random local minimum different from the two homogenous ground states.  
\end{abstract}

{\mbox{}
\vspace{-1.25cm}
\maketitle
}
\vspace{-0.9cm} 

\section{Introduction}
The \emph{disordered Curie--Weiss (CW) model} is a mean-field random ferromagnet defined as follows: consider the complete graph on $N$ vertices and for every edge $(i,j)$, assign a random coupling value $J_{ij}=J_{ji}$ i.i.d.\ according to some non-negative distribution $\mu$ (i.e., $\mbox{supp}(\mu) \in [0,\infty)$). Define the Hamiltonian $H(\sigma)$ for $\sigma \in \{\pm 1\}^N$ by 
\begin{align}
H(\sigma)= - \frac 1N \sum _{1\leq i<j\leq N} J_{ij} \sigma_i \sigma_j\,.
\end{align}
If $\mu$ were instead symmetric on $(-\infty,\infty)$, this would be proportional to the Hamiltonian of the Sherrington--Kirkpatrick spin glass; in our setup, where $\mu$ is non-negative, it is a mean-field analogue of a random ferromagnet. Specifically, when $\mu$ is $\mbox{Bernoulli}(p)$, this corresponds to the Ising model on a dense Erd\H{o}s--R\'enyi random graph, referred to as the \emph{dilute Curie--Weiss model}. In the latter case, the thermodynamics of the model behaves similarly to that of the homogenous Curie--Weiss model, but the coupling disorder induces a more complex energy landscape, whose influence may only manifest itself in the zero-temperature dynamics. The goal of the present paper is to more precisely understand these effects (see Theorem~\ref{mainthm:1}) and relate them to the random constraint satisfaction problem of finding local MINCUTs of a random graph (see \S\ref{sec:local-mincut}).  

The zero-temperature limit of the Glauber dynamics~\cite{Gl63} is called the \emph{zero-temperature dynamics} and in discrete time, is the Markov chain $(X_t)_{t\in \mathbb N}$ with transition matrix given by the following: consider a configuration $\sigma\in \{\pm 1\}^N$; for every $i$, if $\sigma^{(i)}$ is the configuration with $\sigma^{(i)}_j=\sigma_j$ for all $j\neq i$ and $\sigma^{(i)}_i=-\sigma_i$, then
\begin{align}
 P(\sigma,\sigma^{(i)})= \frac 1N \left [\boldsymbol 1\{H(\sigma^{(i)})<H(\sigma)\}+\frac 12 \boldsymbol 1\{H(\sigma^{(i)})=H(\sigma)\}\right ]\,.
 \end{align}
The zero-temperature dynamics is a random walk on the hypercube $\{\pm 1\}^N$, that only assigns positive transition rates to moves that do not increase the energy. If the zero-temperature dynamics stops (i.e., is absorbed) in some particular state, the only possible such \emph{absorbing states} are \emph{local minima} and uniform ground states of $H$. For any non-negative (ferromagnetic) coupling distribution $\mu$, the uniform ground states of the disordered CW model are $\sigma=(1,...,1)$ and $\sigma=(-1,...,-1)$. We call a state $\sigma\in \{\pm 1\}^N$ a \emph{local minimum} of $H$ if it is not a uniform ground state $\pm (1,...,1)$ and is such that for any $\sigma'\in \{\pm1\}^N$ with Hamming distance $d(\sigma,\sigma')=1$, it satisfies $H(\sigma')\geq H(\sigma)$. Notice that if $\mu$ has atoms, there could be connected sets of local minima with the same energy, amongst which the zero-temperature dynamics jumps for all sufficiently large time, in which case we call the set an absorbing set of configurations.

Henceforth, we mostly restrict our attention to the dilute Curie--Weiss model, whose analysis ends up being more straightforward than the general $J_{ij} \geq 0$ case, by letting $J_{ij}\sim \mbox{Ber}(p)$ for $p\in (0,1)$---the case $p=1$ corresponds to the classical Ising Curie--Weiss model. The Ising model on sparse random graphs ($p=p_N\ll \log N/N$) is an extensively studied model with rich relations to random optimization and the theory of spin glasses (e.g.~\cite{DeMo10a} and for a more extensive overview,~\cite{DeMo10}); in the context of zero-temperature dynamics, when the random graph has disjoint components, there are many non-trivial absorbing states, and in particular, the zero-temperature dynamics will have positive probability of absorbing into a local minimum (see e.g.,~\cite{Hag} for a rigorous analysis of zero-temperature dynamics for the Ising model on sparse random graphs). On dense random graphs, the thermodynamics of the Ising model is essentially the same as that of the classical (homogenous) Curie--Weiss model. However, the zero-temperature dynamics is particularly sensitive to small changes in the energy landscape as it can absorb in any local minimum it encounters.

Let $\mathbb P_{\sigma(0)}$ be the product measure over initial configurations, $\mathbb P_{\mathcal J}$ be the product measure over couplings, and $\mathbb P_{\omega}$ be the distribution over the evolution of the Markov chain, with corresponding expectations, $\mathbb E_{\sigma(0)}$, $\mathbb E_{\mathcal J}$, and $\mathbb E_{\omega}$; we will also sometimes write $\mathbb P_{\mathcal J,\omega}=\mathbb P_{\mathcal J} \otimes \mathbb P_{\omega}$ etc. In the present paper we show that while we believe there exist local minima in the energy landscape of the disordered CW model with general ferromagnetic $\mu$ (so that the disorder does indeed make the energy landscape nontrivial---see Question~\ref{question:1}), in the case of the dilute CW model, the zero-temperature dynamics avoids them with high probability, i.e., with probability going to $1$ as $N\to\infty$ (w.h.p.).

\begin{maintheorem}\label{mainthm:1}
For $\epsilon>0$ sufficiently small, the zero-temperature dynamics $(X_t)_{t\geq 0}$ of the Ising model on the Erd\H{o}s--R\'enyi random graph $\mathcal G(N,p)$ with $p>0$ fixed satisfies the following: for every initial configuration $X_0$ with magnetization $\sum_i X_0(i) \geq N^{\frac 12-\epsilon}$, 
\begin{align}
\lim_{N\to\infty} \mathbb P_{\mathcal J,\omega}\left (\lim_{t\to\infty} X_t=(1,...,1) \right)=1\,.
\end{align}
\end{maintheorem}

Note that by heuristic considerations, one would suspect the same should hold for any non-negative coupling distribution $\mu$ that has, for instance, all exponential moments finite. For such distributions, it seems the main obstacle in adapting our argument to give the same result is obtaining some conditional concentration for sums of couplings via an analogue of Proposition~\ref{prop:condition-all-satisfactions}, though we often make use of the convenience of dealing with bounded random variables.

In contrast, if $\{J_{ij}\}_{ij}$ are i.i.d.\ heavy-tailed random variables we prove that with probability bounded away from $0$, the dynamics gets stuck in the exponentially many non-trivial local minima. Results for heavy-tailed couplings are found in~\S\ref{sec:heavy-tail}, Theorem~\ref{thm:local-minima}.

In order to study the zero-temperature dynamics from a physical point of view, a \emph{zero-temperature dynamical order parameter} was introduced~\cite{NNS} and has been examined numerically in a number of models with and without disorder~\cite{NSoverview,NNS,GMNSY17}. Let $\sigma_i(t)$ be the $i$'th spin value of $X_t$, i.e.\ $X_t(i)$. The quantity
\begin{equation}\label{eq:dynamical-order-parameter}
q_D(N)= \mathbb E_{\sigma(0),\mathcal J} [(\mathbb E_\omega[\sigma_1(\infty)])^2]
\end{equation} is designed to capture how much the absorbing state depends on the initial state and how much it depends on the realized evolution of the dynamics. In the above, and throughout the paper, $\sigma_i(\infty)$ denotes $\lim_{t\to \infty} \sigma_i(t)$ if it exists and similarly with $X_\infty(i)$. 

Another way to view the dynamical order parameter $q_D(N)$ is to fix a coupling configuration, and consider a replicated dynamics wherein two replicas $\sigma(t),\sigma'(t)$ undergo zero-temperature dynamics independently ($\omega,\omega'$) from the same initial configuration; one can then ask about the average (under $\mathbb P_{\omega,\omega'}$) evolution of their overlap. The dynamical order parameter is the expectation (under $\mathbb P_{\sigma(0),\mathcal J}$) of this average overlap:
\begin{align}
\mathbb E_{\sigma(0),\mathcal J} \Big[\mathbb E_{\omega,\omega'} \big[N^{-1} \sum _{j=1}^N \sigma_j(\infty) \sigma_j'(\infty)\big]\Big] = q_D(N)\,.
\end{align}

\begin{corollary}\label{cor:qd}
In the dilute CW model, the dynamical order parameter~\eqref{eq:dynamical-order-parameter} has 
\begin{align}
\lim_{N\to\infty} q_{D}(N)=1\,.
\end{align}
In fact, we have the stronger result that,
\begin{align}
\lim_{N\to\infty} \mathbb E_{\sigma(0)}\left[ (\mathbb E_{\mathcal J,\omega}(\sigma_1(\infty)))^2\right]=1\,.
\end{align}

\end{corollary}

We briefly mention that in the above literature regarding this dynamical order parameter, different graphs have also been considered. In~\cite{GMNSY17}, the zero-temperature dynamics of the disordered Curie--Weiss model and random ferromagnet on $\mathbb Z^d$ were studied numerically and heuristically (with light-tailed coupling distributions). There, it was predicted that the zero-temperature dynamics of the disordered Ising model on $(\mathbb Z/N \mathbb Z)^d$ starting from a single state ends up randomly in one of many almost-orthogonal local minima in whose basins of attraction it lies. The numerics suggest that on $d$-dimensional torii, $\lim_{d\to\infty} \lim_{N\to\infty} q_D (N)=0$; combined with Corollary~\ref{cor:qd}, this suggests the existence of a singularity in the $d\to\infty$ behavior of $q_D(N)$. We also note that in the physics literature there has been recent interest in dynamics of the Ising model on networks at low and zero temperatures (e.g.,~\cite{DaSe05} study, at the level of physics, absorption and persistence in a densely-connected small world network).

Over the last several years, significant progress has been made in understanding random optimization problems like extremal cuts on random graphs and random instances of constraint satisfaction problems (see e.g.,~\cite{DMS15,DSS15}) via both heuristic and rigorous connections to spin glasses and other models with quenched disorder. We now discuss a different perspective, related to such random optimization problems, on the fundamental questions underlying Theorem~\ref{mainthm:1}, which may be of independent interest.

\subsection{The local MINCUT problem}\label{sec:local-mincut}
 Consider a dense Erd\H{o}s--R\'enyi random graph $G\sim \cG(N,p)$ for $p>0$ fixed; for any subset $A\subset \{1,...,N\}$, define $\mbox{CUT}_G(A)$ as the number of edges between $A$ and $A^c$. A \emph{local MINCUT} is a partition $(A,A^c)$ of $\{1,...,N\}$ such that for every $A'$ which consists of the addition or removal of one vertex to or from $A$ (in which case we say the Hamming distance $d(A,A')=1$), $\mbox{CUT}_G(A') \geq \mbox{CUT}_G (A)$. A \emph{nontrivial local MINCUT} is one in which both $A$ and $A^c$ are nonempty.  

\begin{question}\label{question:1}
Consider $G\sim \mathcal G(N,p)$ with $p>0$ fixed. Is it the case that 
\begin{align}
\lim_{N\to\infty} \mathbb P(\exists A:A\neq \emptyset, A^c \neq \emptyset, \mbox{CUT}_G(A)\mbox{ is a local MINCUT})=1?
\end{align}
\end{question}

Note that since the Erd\H{o}s--R\'enyi graph is dense, all vertices are connected and have high degrees of dependence; if we were considering a graph with multiple connected components for instance it is obvious how to construct nontrivial local MINCUTs. The requirement that a set $A$ be a local MINCUT is equivalent to demanding that every $v\in A$ have more edges to vertices in $A$ than in $A^c$, and similarly for vertices in $A^c$.

Viewed from this perspective, Theorem~\ref{mainthm:1} can be restated in terms of a greedy search for a local MINCUT defined as follows: start from a uniformly randomly chosen partition $(A,A^c)$ and at every iteration, select a vertex uniformly at random and move it either to $A$ or $A^c$ depending on which move has the lower $\mbox{CUT}_G$ value (if the cut-value is unchanged flip a coin to determine whether to move it). 
\begin{corollary}\label{cor:local-mincut}
With probability going to $1$ as $N\to\infty$, the greedy search for a local MINCUT of $G\sim \cG(N,p)$ with $p>0$ fixed, started from a uniformly random partition of $\{1,...,N\}$ terminates in the trivial partition $(\emptyset, \{1,...,N\})$.\end{corollary}

This suggests the interesting situation where there exist nontrivial local MINCUTs while, with high probability, they are not found by a greedy search algorithm. In fact, there is some numerical evidence in this direction to appear in~\cite{GNSW}. We call such metastable states, \emph{invisible local minima} as they typically do not affect the natural dynamics, even at zero-temperature.

Of course, as with the dilute CW model, the local MINCUT problem can be presented in greater generality by assigning edges of the complete graph i.i.d.\ random weights $w_{ij}$ and asking the analogous questions about $\mbox{CUT}_w (A) = \sum_{i\in A,j\notin A} w_{ij}$. In the case where $w_{ij}$ are symmetric, this corresponds to finding local energy minima of the canonical Sherrington--Kirkpatrick spin glass.  This is an extensively studied question, both at the physics level~\cite{TaEd} and more rigorously recently in the related problem of multiple peaks~\cite{Chat,DEZ}: there the energy landscape is expected to be complex with exponentially many local minima in the system size---~\cite{ABA13, ABC13} developed a complete understanding of the critical points and complexity of the energy landscape in the case where the state space is relaxed to the sphere in dimension $N$. In the above cases, the rugged energy landscape arises due to frustration, a phenomenon that does not exist in the ferromagnetic setup. We also note that in a similar setup to ours, the algorithmic complexity of the local MAXCUT problem has been studied (most recently in~\cite{ABPW}), though that problem again has a very different flavor due to the absence of dominant trivial ground states.

If we instead restrict ourselves to ferromagnetic disorder ($w_{ij}$ are a.s.\ non-negative), we expect that Corollary~\ref{cor:local-mincut} and the techniques of this paper extend to the general case when $w_{ij}$ have light (e.g., exponential or Gaussian) tails; there is again some numerical evidence in this direction~\cite{GNSW}. In contrast, our results in \S\ref{sec:heavy-tail} on heavy-tailed disorder imply an affirmative answer to Question~\ref{question:1} while showing that there, with strictly positive probability, the greedy search terminates in a non-trivial local MINCUT. 

\subsection{Notation}
We introduce some notation that we will use throughout the paper. We say two sequences $f_N$ and $g_N$ are such that  $f_N\lesssim g_N$ if there exists $C>0$ such that $f_N \leq Cg_N$ for all $N$ and we say $f_N \asymp g_N$ if $f_N\lesssim g_N\lesssim f_N$. Finally, we write $f_N = O(g_N)$ if $f_N\lesssim g_N$ and $f_N = o(g_N)$ if $f_N \leq c g_N$ for large enough $N$ for every $c>0$. For readability, we will, throughout the paper, omit floors and ceilings, though all our variables will be integer-valued. We will also assume $N$ is sufficiently large.

The discrete-time zero-temperature dynamics chain is alternately denoted by $(X_t)$ and $\sigma(t)= (\sigma_1(t),...,\sigma_N(t))$ where $t$ is always understood to be integer---clearly, the results of Theorem~\ref{mainthm:1} and Corollary~\ref{cor:qd} would also hold for the analogously defined continuous-time zero-temperature dynamics. The magnetization at time $t$ is given by
\begin{align}
M_t =\sum _{i=1}^N \sigma_i(t)\,.
\end{align}
The effective field on site $i$ at time $t$ is given by
\begin{align}
m_i(t)= \sum _{j\neq i, j=1,...,N} J_{ij} \sigma_j(t)\,.
\end{align}
It will be notationally useful to define the related $\bar m_i(t)=\mbox{sgn}(m_i(t))$ so that $\bar m_i(t)\in\{\pm 1,0\}$ where $\bar m_i(t)=0$ if $m_i(t)=0$. 
Then we let $\{\mathcal S_i(t)\}_{i=1}^\infty$ be the set of \emph{satisfaction random variables}
\begin{align}
S_i(t)=\bar m_i(t)\sigma_i(t)\,,
\end{align}
so that when the dynamics chooses a site $i$ to update, $\sigma_i(t)=\sigma_i(t-1)$ with probability $1$ if $S_i(t)=S_i(t-1)=1$, probability $1/2$ if $S_i(t)=0$, and probability $0$ if $S_i(t)=-1$.

We will use the probability measure $\mathbb P$ to denote the product measure $\mathbb P_{\mathcal J} \otimes \mathbb P_{\omega}$, since we will always be fixing $X_0=\sigma(0)$ and sometimes averging over $\mathcal J,\omega$ at once.

\subsection{Proof approach}
Here we give an overview of our approach to proving Theorem~\ref{mainthm:1}. To avoid the difficulties present in analyzing systems with quenched disorder, particularly with non-Gaussian disorder, our analysis of the zero-temperature dynamics reveals only partial information about the couplings as the dynamics proceeds. This gradual ``revealing scheme" may be of independent interest in analyzing the short-time ($t\ll N$) dynamics of other systems with quenched disorder. Here it allows us to bound the drift of the magnetization chain $(M_t)_{t\in \mathbb N}$ conditioned on this partial information from below, and compare $(M_t)_{t\in \mathbb N}$ to a random walk with positive drift. More precisely,

\begin{itemize}
\item In \S 2.1, we define a \emph{revealing scheme} to see the evolution of the zero-temperature Markov chain as measurable w.r.t.\ the $\sigma$-algebra $\mathcal F_t$ generated by the sequence $(Y_t)_t$ of updated sites, their satisfactions at update time $(S_{Y_t}(t))_t$, and the couplings $(\{J_{Y_t,Y_k}\}_{k=t+1,...,T})_t$ for $T=N^{\frac 12+\delta}$. To gradually reveal this information, we first fix the update sequence $\{Y_1,...,Y_T\}$ then sequentially jointly reveal the satisfaction and aforementioned couplings of the next site to update.

\item In \S 2.2, we compare the joint distribution of $(\{J_{Y_t,Y_k}\}_{k=t+1,...,T})_t$ given $S_{Y_t}(t)$ and $\cF_{t-1}$ to a product measure. In particular, for all short times $t=O( N^{\frac 12+\delta})$ we show that the conditional joint law of the couplings revealed at time $t$ dominates i.i.d.\ $\mbox{Ber}(p-O(N^{\frac 12+2\delta}))$ and is dominated by i.i.d.\ $\mbox{Ber}(p+O(N^{\frac 12+2\delta}))$.

\item In \S 2.3, we show that for all times $t\leq N^{\frac 12+\delta}$, as long as $M_t \geq N^{\frac 12-\delta}$ holds, the chain $(M_t)_{t\geq 0}$ has a positive drift of at least $cN^{-\delta/2}$, so that, $(M_t)_{t\geq 0}$ typically stochastically dominates a random walk with drift of $c N^{-\delta/2}$. 

\item By this comparison, at time $T=N^{\frac 12+\delta}$, $M_T \geq N^{\frac 12+\frac \delta 2}$ and all sites that have not yet been updated have a positive effective field. This can then be boosted in \S 2.4 to show that by time $T'=N^{2/3}$ all sites have a positive field and the dynamics will thereafter quickly absorb into the all-plus ground state.
\end{itemize}

\section{The Dilute Curie--Weiss model}

\subsection{Random mapping representation}
Denote by $Y_k\in \{1,...,N\}$ the site chosen at time step $k$ to be updated so that for every $k\in \mathbb N$, $\mbox{dist}(Y_k)=\mbox{Uni}(\{1,...,N\})$. Then consider the sequence of update sites and their satisfactions, $(Y_k, S_{Y_k}(k))_{k}$.

In the dilute Curie--Weiss model, $\mu$ is atomic and there may be zero-energy flips with positive probability, so $(X_t)_{t\geq 0}$ is not measurable with respect to the sigma-algebra generated only by $X_0$ and $(Y_k,S_{Y_k}(k))_{k\leq t}$. Thus define a sequence of i.i.d.\ random variables $(B_k)_{k\geq 0}$ with $B_1\sim 2\mbox {Ber}(1/2)-1$ which will determine the spin at site $Y_k$ in the case that $m_{Y_k}(k)=0$. Then the history of the chain $(X_k)_{k\leq t}$ is fully determined by $X_0=\sigma(0)$ and the sequence $(Y_k,\mathcal S_{Y_k}(k),B_k)_{k\leq t}$. 

The \emph{grand coupling} is the coupling of two independent dynamics with different $X_0$ such that both dynamical realizations use the same random variable sequence $(Y_k)_{k\geq 1}$ and $(B_k)_{k\geq 0}$ as well as the same couplings $\{J_{ij}\}$. The grand coupling has the added feature that it preserves monotonicity, so that if $X_0\succeq X_0'$ then $X_k\succeq X_k'$ for all $t\geq 0$. Moreover, by permutation invariance of the measure $\mathbb P_{\cJ,\omega}$, we can begin by permuting all initial configurations so that if $M_0 \geq M_0'$, $X_0 \succeq X_0'$ and identify initial configurations only with their magnetization, then apply the grand coupling of $\cJ,\omega$.

\subsection{Preliminary estimates}

We first estimate the probability of the magnetization of $X_0$ being atypical. The following is a consequence of e.g., Berry--Esseen theorem.

\begin{fact}\label{fact:init-config} For every $\epsilon>0$, we have
\begin{align}
\mathbb P_{\sigma(0)} (X_0:|M_0|\geq N^{\frac 12-\epsilon})=1-O(N^{-\epsilon})\,.
\end{align}
\end{fact}

By the grand coupling and permutation symmetry of the model, it suffices to prove Theorem~\ref{mainthm:1} for a fixed $X_0$ with $M_0=N^{\frac 12-\epsilon}$; for every coupling and dynamical realization in which the chain with that initial configuration absorbs in $X_\infty =(1,...,1)$, via the grand coupling, so does every chain with more positive initial configuration. Moreover, in order to simplify our considerations, by monotonicity we may assume that for all $k\geq 1$, we always have $B_k=-1$; via the grand coupling of the dynamics, it suffices to show that this chain absorbs into the all-plus configuration w.h.p.\ to prove Theorem~\ref{mainthm:1}. Abusing notation, $(X_k)_k$ will henceforth refer to this new chain with $B_k=-1$ for all $k$. In conjunction with this change, whenever $m_{Y_k}(k)=0$, we set $S_{Y_k}(k)=-\sigma_{Y_k}(k-1)$.

\begin{definition}[Revealing scheme]\label{def:revealing-scheme}
Fix $T=N^{\frac 12 +3\epsilon}$ and for every $t\in \mathbb N$, let $\{\mathcal F_t\}_t$ be the filtration of $\sigma$-algebras generated by 
\begin{align}
(Y_1,...,Y_{T\vee t}), (S_{Y_k}(k))_{k\leq t}, \mbox{ and } (J_{Y_k Y_l})_{k \leq t, l\leq T}
\end{align}
The chain $X_t$ is measurable w.r.t.\ the $\sigma$-algebra generated by $(Y_k, S_{Y_k}(k))_{k=1,...,t}$ and $X_0$ (and therefore measurable w.r.t.\ $\mathcal F_t$ and $X_0$).  
\end{definition}

It will be crucial to understand the conditional distribution of $J_{ij}$ given $\mathcal F_t$ when $i\in \{Y_k\}_{k\leq t}$ but $j\notin \{Y_k\}_{k\leq t}$, as well as the joint law of such $\{J_{ij}\}_j$. (This is the main reason we restrict ourselves to Bernoulli $J_{ij}$, where the distribution is determined by the mean, instead of general light-tailed random variables like half-normal distributed random variables. To extend Theorem~\ref{mainthm:1} to that setting, the main technical hurdle is obtaining appropriate analogues to Proposition~\ref{prop:condition-all-satisfactions} for exponential moments in order to obtain conditional concentration of sums of couplings, namely~\eqref{eq:coupling-concentration}.)

We will need the following notation: for any sequence of order updates, denote by $\mathcal R_t$ the set of vertices whose clocks have rung more than once before time $t$, i.e., 
\begin{align}\label{eq:R-t}
\mathcal R_t= \Big\{i:\sum_{k\leq t} \boldsymbol 1\{Y_k=i\}>1\Big\}\,.
\end{align}
\begin{proposition}\label{prop:condition-all-satisfactions} Let $t \leq T=N^{\frac 12+3\epsilon}$, suppose that $j=Y_t$, $j\notin \mathcal R_t$, that $\ell \notin \{Y_k\}_{k=1}^t$, and suppose that $|M_0| \leq N^{\frac 12+3\epsilon}$. Then we have that, 
\begin{align}\label{eq:corr-want-to-show}
\left |\mathbb E_{\cJ} \left[J_{\ell j} \given \mathcal F_{t-1}, S_{j}(t), \{J_{ja}\}_{a\in \{Y_k\}_{t+1}^T - \mathcal R_T-\{\ell\}}\right] -p\right| = O(N^{-\frac 12+4\epsilon})\,.
\end{align}
\end{proposition}
\begin{proof}
First fix the update sequence $\{Y_k\}_{k=1}^T$ and fix any such $t$ and $\ell$. 
Let 
\begin{align}
\tilde{\mathbb E} [\,\cdot\,] & = \mathbb E_{\cJ}[\,\cdot \mid \mathcal F_{t-1}, \{J_{ja}\}_{a\in \{Y_k\}_{k=t+1}^T-\mathcal R_T-\{\ell\}}]\,.
\end{align}
Now we can expand 
\begin{align}\label{eq:tilde-expansion}
\tilde {\mathbb E}[m_{j}  (t)  \mid S_{j}(t)]-  \mathbb E_{\cJ} [m_{j}(0)]= &  \,\,\, \sigma_{\ell}(0) \tilde { \mathbb E}[J_{\ell j}-p \mid S_j(t)] \\ 
& + \sum_{i \in \{Y_k\}_{k=1}^{T}-\{\ell\}} \big[\sigma_{i}(t)  \tilde {\mathbb E}[J_{i j}\mid S_j(t)] -p\sigma_i(0)\big]  \nonumber \\
  &  + \sum_{i\notin \{Y_k\}_{k=1}^T\cup \{\ell\}} \sigma_i(0)\tilde {\mathbb E}[J_{i j}-p \mid S_j(t)] + O(|\cR_T|)\,.\nonumber
\end{align}
Since $|\mathcal R_T| \leq T$, we can replace $O(\cR_T)$ by $O(T)$; moreover the second sum consists of at most $T$ terms and is thus bounded in absolute value by $2T$. Now suppose without loss of generality that $m_j(t) > 0$ and $\sigma_{\ell} (0)=+1$---the same argument carries through in the other cases. We claim, first of all that for every $m\notin \{Y_k\}_{k=1}^T$ with $\sigma_m(t)=-1$,
\begin{align}\label{eq:conditioning-decrease}
\tilde {\mathbb E} [ J_{m j}\mid m_{j}(t)>0] \leq \tilde {\mathbb E} [ J_{m j}]=p\,.
\end{align}
To see this, we write by Bayes' Theorem, the left hand side above as 
\begin{align}
\tilde {\mathbb P}(m_{j}(t)>0 \mid J_{m j}=1) \tilde {\mathbb P} (J_{m j} = 1) \tilde {\mathbb P} (m_{j}(t)>0)^{-1}\,.
\end{align}
Observe that since $\tilde {\mathbb P}$ doesn't condition on $m_{j}$ nor on $m_{m}$ at any time, $J_{m j}$ is independent of the $\sigma$-algebra conditioned on under $\tilde {\mathbb P}$, and thus its conditional distribution is $\mbox{Ber}(p)$ so that $\tilde {\mathbb E}[J_{m j}]=p$. Moreover, expanding out $m_{j}(t)>0$, we see that under $\tilde {\mathbb P}$, all the summands except $\sigma_m(0)J_{mj}$ are conditionally independent of $J_{mj}$, so that because $\sigma_m(t)=-1$, conditioning also on $J_{mj}=1$ only decreases $m_{j}(t)$, implying~\eqref{eq:conditioning-decrease}. Analogously, if $\sigma_{m}(0)=1$, and $m\notin \{Y_k\}_{k=1}^T$, then $\tilde{\mathbb E}[J_{mj} \mid m_{j}(t)>0] \geq p$. Therefore, every summand in the third term in the right-hand side of~\eqref{eq:tilde-expansion} is nonnegative. 
As a result, we have
\begin{align}\label{eq:positivity}
0 \leq \sigma_\ell(0) \tilde {\mathbb E}[J_{\ell j} -p\mid S_j(t)]+ \sum_{i\notin \{Y_k\}_{k=1}^T} & \, \sigma_i(0) \tilde{\mathbb E}[J_{ij} - p\mid S_j(t)] \nonumber \\
& \leq |\tilde{\mathbb E} [m_j(t) \mid S_j(t)]|+ |\mathbb E[m_j(0)]|+O(T)\,.
\end{align}

We now upper bound the right-hand side of~\eqref{eq:positivity}. Writing each term out, we see that under $\tilde {\mathbb P}$, except for at most $t=O(N^{\frac 12+3\epsilon})$ summands, $m_j(t)$ is distributed as a difference of two binomial random variables with mean that is $M_0 \pm O(N^{\frac 12+3\epsilon})$ and variance $O(N)$. (Under the measure $\mathbb P_{\cJ}$, $m_j(0)$ is of course just a difference of two binomials.) In that case, standard lower and upper tail estimates, via Chernoff bounds for binomial random variables, along with the fact that $|M_0|\leq N^{\frac 12 +3\epsilon}$ imply that 
\begin{align}
\tilde {\mathbb E} [m_{j}(t) \mid S_j(t)]  & =O(N^{\frac 12 +4\epsilon})\,, \qquad \mbox{and} \qquad  {\mathbb E}_{\cJ}[m_{j}(0)]  =O(N^{\frac 12 +4\epsilon})\,.
\end{align}

Returning to~\eqref{eq:positivity}, by~\eqref{eq:conditioning-decrease} we see that when $\sigma_\ell(0)=+1$, we have
\begin{align}
\tilde{\mathbb E}[J_{\ell j} - p\mid S_j(t)] + \sum_{i\notin \{Y_k\}_{k=1}^T, \sigma_i(0) = +1} \tilde {\mathbb E}[J_{ij} - p\mid S_j(t)] = O(N^{\frac 12+4\epsilon})\,.
\end{align}
There are deterministically at least $N/3$ terms in the sum above, so that in order to conclude, it suffices to show that for an $i\notin \{Y_k\}_{k=1}^T$ with $\sigma_i(0)=1$, we have that $\tilde {\mathbb E} [J_{\ell j} \mid S_j(t)]= \tilde{\mathbb E}[J_{ij} \mid S_j(t)]$.
(One would then divide both sides by the number of terms in the sum above, and obtain the desired.) If $\ell\notin \{Y_k\}_{k=1}^T$, this is evident by symmetry. Else, write, 
\begin{align}
\tilde{\mathbb E} [J_{\ell j} \mid m_j(t)>0] = \frac{\tilde{\mathbb P}(J_{\ell j}=1)\tilde{\mathbb P}(m_j(t)>0\mid J_{\ell j}=1)}{\tilde{\mathbb P}(m_j(t)>0)}\,.
\end{align}
The conditioning in $\tilde {\mathbb P}$ is independent of $J_{\ell j}$ as well as $J_{ij}$ so the first term in the numerator is the same under $\ell \mapsto i$. Likewise, for the other two terms, expanding out $m_j(t)$, one sees that since $\tilde {\mathbb P}$ conditions on $J_{kj}$ for every $k \in \{Y_1,...,Y_{t-1}\}$, the other couplings conditioned on under $\tilde{\mathbb P}$ do not affect the distribution of $m_j(t)$; moreover $\sigma_{i}(t)=\sigma_{\ell}(t)=\sigma_i(0)=\sigma_\ell(0)$, so that those remaining two terms are also unchanged under $\ell\to i$. 
\end{proof}

\subsection{Short time dynamics}\label{sec:short-time}
In this section we study the evolution of $(X_t)_{t\in \mathbb N}$ started from $X_0$ such that $M_0 =  N^{\frac 12 -\epsilon} $ until $T= N^{\frac 12 +3\epsilon} $. 

We begin by fixing the update sequence $Y_1,...,Y_T$, then for each realization of the update sequence, we bound probabilities of evolutions of the chain under $\mathbb P_{\mathcal J}$. 
Recall the definition of $\mathcal R_k$ from~\eqref{eq:R-t}. We define the following good events for the update sequence $\{Y_k\}_{k=1}^T$:
\begin{align}\label{eq:two-conditions}
\Gamma_{\omega,t}^1:=\left\{|\mathcal R_{t}| < N^\epsilon \vee \frac {2t^2}{N}\right\}\qquad \mbox{and}\qquad \Gamma_{\omega,t}^2:=\bigg\{|\sum_{k=1}^{t} \sigma_{Y_k} (0)|< t^{\frac 12+\epsilon}\bigg\}\,.
\end{align}
The dynamical good event is then defined as $\Gamma_\omega= \bigcap_{t=N^{\epsilon}}^{T} \Gamma_{\omega,t}^1 \cap \Gamma_{\omega,t}^2$. 

\begin{lemma}\label{lem:good-update-sequence}
Let $X_0$ be such that $M_0 = N^{\frac 12 -\epsilon}$. There exists $c(\epsilon)>0$ such that 
\begin{align}
\mathbb P_{\omega} ( \Gamma_\omega^c) \lesssim Te^{-cN^{2\epsilon^2}}\,.
\end{align}
\end{lemma}

\begin{proof}
We union bound over the $(\Gamma_{\omega,t}^i)^c$ for $i=1,2$ and $N^{\epsilon}\leq t\leq T$. The bound for $i=1$ follows from a union bound over all $t\leq T$ and the following. For each $t\leq T$, the probability of selecting a site that has already been updated at time $t$ is at most $\frac tN \leq N^{-\frac 12 +3\epsilon}$.  Therefore, $|\mathcal R_t| \preceq \mbox{Bin} (t, t/N)$ so that by the Chernoff inequality, 
\begin{align}
\mathbb P_\omega\big(|\mathcal R_t|\geq \max\{N^{\epsilon}, 2t^2 N^{-1}\}\big) \leq 2e^{-cN^{\epsilon}}\,,
\end{align}
for some $c>0$.
For the $i=2$ bound, because the update order $(Y_k)_{k\leq T}$ is independent of $\sigma(0)$ and $t\leq N^{\frac 12 +3\epsilon}$, by Hoeffding's inequality,
\begin{align}
\mathbb P_\omega \bigg(|\sum_{k=1}^{t} \sigma_{Y_k} (0)|\geq t^{\frac 12+\epsilon}\bigg)\leq 2e^{-ct^{2\epsilon}}\,,
\end{align}
because uniformly in $t\leq T$, the probabilities of $\sigma_{Y_t}(0)=\pm 1$  are within $\frac {M_0}{N}$ of $\frac 12$ and $M_0= N^{\frac 12-\epsilon}$. Union bounding over all $N^{\epsilon}\leq t\leq T$ yields the desired. 
\end{proof}

By Lemma~\ref{lem:good-update-sequence}, without loss, we can now restrict our attention to realizations of $\omega$ such that the sequence $Y_1,...,Y_T$ satisfies $\Gamma_{\omega}$.

We also define a good coupling event measurable with respect to $\mathcal F_{t}$ which we will restrict our attention to. Fix a sequence $(Y_1,...,Y_T)$. Begin by defining
\begin{align}\label{eq:def-z}
Z_{++}(t)= \big\{j\in \{Y_k\}_{k=1}^{t-1} - \mathcal R_T\,:\, \sigma_j (0)=+ 1, \sigma_j(t)=+1\big\}
\end{align}
and the analogously defined $Z_{--},Z_{+-},Z_{-+}$, measurable with respect to $(Y_k,S_{Y_k}(k))_{k< t}$. Then for $N^\epsilon \leq t \leq T$, let
\begin{align}
\Gamma_{\mathcal J,t}^{++} = \bigcap_{s=N^{\epsilon}}^{t} \Big\{ \Big|\sum_{j\in Z_{++}(s)} [J_{Y_s j} -p]\Big| \leq s^{\frac 12+\epsilon}  \Big\}\,.
\end{align} 
Then let $\Gamma_{\mathcal J,t}= \Gamma^{++}_{\mathcal J,t} \cap \Gamma^{+-}_{\mathcal J,t} \cap \Gamma^{-+}_{\mathcal J,t} \cap \Gamma^{--}_{\mathcal J,t}$.

\begin{lemma}\label{lem:coupling-concentration}
Fix an initial configuration $X_0$ such that $M_0 = N^{\frac 12 - \epsilon}$. There exists $c(p)>0$ such that for every sequence $Y_1,...,Y_T$, and every  $N^{\epsilon} \leq t\leq T$,
\begin{align}
\mathbb P_{\mathcal J} ( \Gamma^c_{\mathcal J,t}) \lesssim te^{-cN^{2\epsilon^2}}\,.
\end{align}
\end{lemma}
\begin{proof}
By a union bound over all $N^{\epsilon} \leq s \leq t$, it suffices to show that there exists $c(p)>0$ such that for every such $s$,
\begin{align}\label{eq:coupling-concentration}
\mathbb P_{\mathcal J} \Big( \Big| \sum_{j\in Z_{++}(s)} [J_{{Y_s} j}-p]\Big| \geq s^{\frac 12 +\epsilon} \given (Y_k,S_{Y_k}(k))_{k<s}\Big) & \lesssim \exp(-cs^{2\epsilon}) 
\end{align}
and similarly for $Z_{+-}(s), Z_{-+}(s), Z_{--}(s)$, then average over all possible realizations of $(S_{Y_k}(k))_{k<s}$. We use Proposition~\ref{prop:condition-all-satisfactions} to jointly sample the couplings $\{J_{Y_s j}\}_{j\in Z_{++}(s)}$ and $S_{Y_1}(1),...,S_{Y_{s-1}}(s-1)$ according to the revealing scheme defined in Definition~\ref{def:revealing-scheme}. One sees by Proposition~\ref{prop:condition-all-satisfactions} that under this revealing process, independently of the couplings that have already been revealed and the information obtained from $(S_{Y_k}(k))_{k<s}$, the distribution of $J_{j Y_s}$ dominates $\mbox{Ber}(p-cN^{-\frac 12+4\epsilon})$ and is dominated by $\mbox{Ber}(p+cN^{-\frac 12+4\epsilon})$.  Therefore, conditional on any actualization of $(S_{Y_k}(k))_{k<s}$ and therefore $Z_{++}(s)$, the joint distribution of $\{J_{Y_s j}\}_{j\in Z_{++}(s)}$ is dominated by a product measure of $\mbox{Ber}(p+cN^{-\frac 12+4\epsilon})$ and dominates a product measure of $\mbox{Ber}(p-cN^{-\frac 12 +4\epsilon})$. At that point, using Chernoff--Hoeffding inequality, we see that~\eqref{eq:coupling-concentration} holds. 
\end{proof}

The main estimate on the effective fields at short times $t\leq T$ is the following.  

\begin{proposition}\label{prop:probability-positive-field}
Fix $X_0$ with magnetization $M_0= N^{\frac 12 -\epsilon}$. Suppose for $t\leq T$, we are on a set in $\mathcal F_{t-1}$ such that $\Gamma_{\omega}\cap \Gamma_{\mathcal J,t}$ holds and $M_{t-1} \geq M_0-N^{2\epsilon}$.
Suppose also that $Y_t \notin \mathcal R_t$. There exists $c(p)>0$ such that for sufficiently small $\epsilon>0$, we have
\begin{align}\label{eq:main-estimate}
 \mathbb P_{\mathcal J} \left(m_{Y_t}(t)>0\mid \mathcal F_{t-1} \right)\geq \frac 12 +cN^{-\epsilon}\,.
\end{align}
In particular, this estimate holds independent of $\sigma_{Y_t}(0)=\sigma_{Y_t}(t-1)$.
\end{proposition}

\begin{proof}
For ease of notation, let $i=Y_t$. 
By definition, the $\sigma$-algebra $\mathcal F_{t-1}$ reveals information about the couplings $\{J_{i Y_k}\}_{k=1}^{t-1}$ and is independent of any other couplings to $i$. 
Using the fact that if $j\notin \{Y_{k}\}_{k<t}$, $\sigma_j(0)=\sigma_j(t)$, the event $\Gamma_{\omega,t}^2$ implies
\begin{align}
\sum_{j\notin \{Y_k\}_{k<t}} \sigma_j(0) = \sum_{j\notin \{Y_k\}_{k<t}} \sigma_j(t)\geq N^{\frac 12-\epsilon}-N^{\frac 14+3\epsilon}\,.
\end{align}

Now expanding $m_i(t)$, we obtain
\begin{align}\label{eq:eff-field}
m_{i}(t)=\sum_{j\in \{Y_k\}_{k< t}} \sigma_j(t) J_{ij}+\sum_{j\notin \{Y_k\}_{k< t}} \sigma_j(0) J_{ij}\,.
\end{align}
We first consider the second sum in~\eqref{eq:eff-field}. As remarked earlier, because $j\notin \{Y_k\}_{k<t}$, conditional on $\mathcal F_{t-1}$, these $J_{ij}$ are distributed as i.i.d.\ $\mbox{Ber}(p)$. Thus, given $\Gamma_{\omega,t}^2$ holds, we have the following conditional (on $\cF_{t-1}$) stochastic domination:
\begin{align}\label{eq:binomials}
\sum_{j\notin \{Y_k\}_{k<t}} \sigma_j(0) J_{ij} &
\succeq \mbox{Bin}\big(\tfrac {N}2+(\tfrac {1-\epsilon}2)N^{\tfrac 12-\epsilon},p\big)- \mbox{Bin}\big(\tfrac {N}2-(\tfrac {1-\epsilon}2)N^{\tfrac 12-\epsilon},p\big)\,.
\end{align}
By Berry--Esseen Theorem, we thus obtain for some $c(p)>0$,
\begin{align}\label{eq:havent-rung-yet}
 \mathbb P_{\cJ}\bigg(\sum_{j\notin \{Y_k\}_{k<t}} \sigma_j(0)J_{ij} \geq \frac p2 N^{\frac 12 -\epsilon}\bigg) \geq \frac 12 +cN^{-\epsilon} + O(N^{-\frac 12})\,.
\end{align}
We now control the contribution of the first term of~\eqref{eq:eff-field}. Consider separately the case $t\leq N^{2\epsilon}$ and $t\geq N^{2\epsilon}$. If $t\leq N^{2\epsilon}$, then clearly the first sum in~\eqref{eq:eff-field} is bounded above by $t \leq N^{2\epsilon}$. Now consider the case $t\geq N^{2\epsilon}$.

Recall the definitions of $Z_{++}, Z_{+-}, Z_{-+}$ and $Z_{--}$ from~\eqref{eq:def-z}. Then we can expand
\begin{align}\label{eq:sum-four-terms}
\sum_{j\in \{Y_k\}_{k<t}} \sigma_j(t) J_{ij}  = & \sum _{j\in \mathcal R_t}\sigma_j(t) J_{ij} + \sum_{j\in Z_{++}(t)} J_{ij} +\sum_{j\in Z_{-+}(t)}J_{ij} \nonumber \\
& - \sum_{j\in Z_{+-}(t)}  J_{ij} -\sum_{j\in Z_{--}(t)}J_{ij}\,.
\end{align}
First of all, by the fact that $\Gamma_{\omega,t}^1$ holds, $|\mathcal R_t|$ and in turn the first sum on the right-hand side of~\eqref{eq:sum-four-terms}, are bounded above by $N^{7\epsilon}$ independently of $t\leq T$.

For the latter four sums in~\eqref{eq:sum-four-terms}, we use the fact that $\Gamma_{\mathcal J,t}$ holds. 
First we bound how many summands there are in each of the four terms. Since $\Gamma_{\omega,t}^1,\Gamma_{\omega,t}^2$ both hold and $M_{t-1}\geq M_0-N^{2\epsilon}$, 
\begin{align}\label{eq:difference-z}
\big||Z_{++}(t)| + |Z_{+-}(t)|-|Z_{--}(t)|-|Z_{-+}(t)|\big| & \leq N^{\frac 14 +3\epsilon}+N^{7\epsilon}\qquad\mbox{and} \nonumber \\
\qquad |Z_{+-}(t)|-|Z_{-+}(t)| & \leq N^{2\epsilon}+N^{7\epsilon}\,,
\end{align}
so that $|Z_{++}(t)|+|Z_{-+}(t)| \geq |Z_{--}(t)|+|Z_{+-}(t)|-N^{\frac 14+3\epsilon}- 2N^{ 7\epsilon}-N^{2\epsilon}$. Now we observe that by definition of $\Gamma_{\mathcal J,t}$, we have 
\begin{align}\label{eq:gamma-j-applied}
 \Big| \sum_{j\in Z_{++}(t)} [J_{{Y_t} j}-p]\Big| \leq t^{\frac 12 +\epsilon} \end{align}
and likewise when summing over $Z_{+-}(t),Z_{-+}(t)$ and $Z_{--}(t)$. 
If $Z_+(t) = Z_{++}(t) + Z_{-+}(t)$ and $Z_-(t) = Z_{--}(t) + Z_{+-}(t)$, by applying~\eqref{eq:gamma-j-applied} separately to $Z_+(t)$ and $Z_-(t)$, we see that under the event $\Gamma_{\omega} \cap \Gamma_{\mathcal J,t}$, by~\eqref{eq:difference-z} and the fact that $N^{2\epsilon} \leq t\leq N^{\frac 12+3\epsilon}$, for $\epsilon$ small,
\begin{align}
\sum_{j\in \{Y_k\}_{k<t}} \sigma_j(t) J_{ij} \geq p[Z_+(t)-Z_-(t)]-4t^{\frac 12+\epsilon} \geq -2N^{\frac 14+3\epsilon}\,,
\end{align}
and if $t\leq N^{2\epsilon}$ the same bound trivially holds.
Combined with~\eqref{eq:havent-rung-yet}, this implies that
\begin{align}
\mathbb P_{\mathcal J}(m_{Y_t}(t)\leq 0 \mid \mathcal F_{t-1})& \leq  \mathbb P_{\mathcal J}\bigg(\sum_{j\notin \{Y_k\}_{k<t}} \sigma_j(0)J_{ij}  \leq 3 N^{\frac 14 +3 \epsilon} \given \mathcal F_{t-1}\bigg) \nonumber \\
& \leq \frac 12 -cN^{-\epsilon}+O(N^{-\frac 12})\,,
\end{align}
yielding the desired for all sufficiently small $\epsilon>0$.
\end{proof}

The following shows the implications of Proposition~\ref{prop:probability-positive-field} for the magnetization chain.
\begin{corollary}\label{cor:effective-field-bound}
Fix $X_0$ with magnetization $M_0=N^{\frac 12-\epsilon}$. Suppose that $t\leq T$, that $\Gamma_{\omega}\cap\Gamma_{\mathcal J,t}$ holds, and that $M_{t-1}\geq M_0 - N^{2\epsilon}$. Suppose also that $Y_t \notin \mathcal R_t$; if $\sigma_{Y_t}(0)=-1$,
\begin{align}
\mathbb P_{\mathcal J} (M_{t}-M_{t-1} = +1\mid \mathcal F_{t-1}) = 1-\mathbb P_{\mathcal J} (M_{t}-M_{t-1} = 0\mid \mathcal F_{t-1}) \geq \frac 12+cN^{-\epsilon}\,,
\end{align} 
for some uniform $c(p)>0$ and if $\sigma_{Y_t}(0)=+1$,
\begin{align}
\mathbb P_{\mathcal J} (M_{t}-M_{t-1} = 0\mid \mathcal F_{t-1}) = 1-\mathbb P_{\mathcal J} (M_{t}-M_{t-1} = -1\mid \mathcal F_{t-1}) \geq \frac 12+cN^{-\epsilon}\,.
\end{align} 
\end{corollary}
 
We now use Corollary~\ref{cor:effective-field-bound} to lower bound the magnetization chain $(M_t)_{t\leq T}$. Define for $\epsilon,t,\theta>0$, the following subset of the filtration $\mathcal F_t$:  
\begin{align}
\mathcal B_{\epsilon,t,\theta}  & =\{ M_k > M_0 - N^{2\epsilon}+\theta k N^{-\epsilon} \mbox{ for all } k\leq t\}\,,
\end{align}
and observe that $\mathcal B_{\epsilon,t,\theta}$ is an  increasing event in the magnetization.

\begin{proposition}[Magnetization lower bound] \label{prop:magnetization-bound}
Fix $X_0$ such that $M_0= N^{\frac 12-\epsilon}$. Suppose $\Gamma_\omega\cap\Gamma_{\mathcal J,T}$ holds for the realization of $\mathcal F_{t}$. Then for all $\theta>0$ sufficiently small, there exists a constant $d(p,\theta)>0$ such that  
\begin{align}\mathbb P_{\mathcal J}\left(\mathcal B^c_{\epsilon,T,\theta} \right)\leq T e^{-dN^{2\epsilon^2}}\,.
\end{align}
\end{proposition}

\begin{proof}
We fix any $t\leq T$ and show $\mathbb P_{\mathcal J} (\cB_{\epsilon,t,\theta}^c) \leq te^{-dN^{2\epsilon^2}}$, implying in particular, the desired. In order to do so, we consider two random walk chains, $\mathcal M^+_t$ which lower bounds the change in magnetization over times when $\sigma_{Y_k}(0)=+1$ and $\mathcal M^-_t$ which does the same over times when $\sigma_{Y_k}(0)=-1$. The chains $\mathcal M^+_t$ and $\cM_t^-$ are defined as follows: let $B_i^+, B_i^-$ be i.i.d.\ $\mbox{Ber}(\frac 12+cN^{-\epsilon})$ for $c(p)>0$ given by Corollary~\ref{cor:effective-field-bound},
\begin{align}
\mathcal M^+_{t+1}-\cM^+_t = 1-B_t^+\,, \qquad \mbox{and}\qquad \cM^-_{t+1}-\cM^-_t = B_t^-\,,
\end{align}
and $\cM_0^+ = \cM_0^-=0$.
Define also $t_+= \sum_{j=1}^t \boldsymbol 1\{\sigma_{Y_j}(0)=+1\}$ and analogously $t_-$. Clearly $t_+ + t_- =t$. As long as $\cB_{\epsilon,{t-1},\theta}$ holds for some $\theta \geq 0$, by Corollary~\ref{cor:effective-field-bound},
\begin{align}\label{eq:rw-domination}
M_t \succeq \mathcal M_t-|\mathcal R_t|\qquad \mbox{where} \qquad \cM_t:= M_0+\mathcal M_{t_+}^+ + \mathcal M_{t_-}^-\,.
\end{align}

Let $\bP$ be the probability distribution over $(B_i^+,B_i^-)_{i}$ and therefore the random walks $\cM_t^+$ and $\cM_t^-$. Moreover, let $\mathcal A_{\epsilon,t,\theta}$ be the analogue of $\mathcal B_{\epsilon,t,\theta}$ for the random walk $\mathcal M_t-|\cR_t|$ (note that $\cR_t$ is fully determined by the sequence $(Y_k)_{k\leq T}$).
We prove that $\bP (\cM_t - |\cR_t| \in \cA_{\epsilon,t,\theta}^c) \leq t e^{-dN^{2\epsilon^2}}$ for sufficiently small $\theta>0$ for some $c(p,\theta)>0$; this would imply the desired since $M_t$ and $\cM_t$ can be coupled through all realizations of $\cA_{\epsilon,t,\theta}$ in a way that $M_t\geq \cM_t$, implying that $\cB_{\epsilon,t,\theta}$ also holds. 

To prove this, note first that for all $t\leq N^{2\epsilon}$, $\cA_{\epsilon,t,\theta}$ holds trivially for every $\theta\geq 0$. Thus consider $t\geq N^{2\epsilon}$. Observe that since $\Gamma_{\omega}$ holds, for every $N^{\epsilon}\leq t\leq T$,
\begin{align}\label{eq:t+t-}
|t_+-t_-| \leq t^{\frac 12+\epsilon}\,,\qquad \mbox{and}\qquad |\cR_t| \leq N^{\epsilon}\vee 2t^2/N\,.
\end{align}
By Hoeffding's inequality, there exists $d>0$ such that for every $\delta>0$,
\begin{align}
\bP \Big(|\cM_{t_+}+t_+(\tfrac 12-cN^{-\epsilon})| & \geq t_+^{\frac 12+\delta}\Big)\leq 2e^{-dt_+^{2\delta}}\,, \\
\bP\Big(|\cM_{t_-}-t_-(\tfrac 12+cN^{-\epsilon})| & \geq t_-^{\frac 12+\delta}\Big) \leq 2e^{-dt_-^{2\delta}}\,.
\end{align}
Combining the above with~\eqref{eq:t+t-}, we see that for some different $d>0$,
\begin{align}
\bP\Big(|\cM_{t_+}+\cM_{t_-} - (\tfrac{t_--t_+}2 + tcN^{-\epsilon})| \geq t_+^{\frac 12+\delta}+t_-^{\frac 12+\delta}\Big) \leq 4e^{-dt^{2\delta}}
\end{align}
and in particular, letting $\epsilon=\delta$, since $|t_+ - t_-| \leq t^{\frac 12+\epsilon}$,
\begin{align}
\bP(\cM_t \leq M_0-3t^{\frac 12+\epsilon}+tcN^{-\epsilon})\leq 2e^{-dt^{2\epsilon}}\,.
\end{align}
Now note that by the bound on $|\mathcal R_t|$ in~\eqref{eq:t+t-}, for all $N^{\epsilon}\leq t\leq T$, 
\begin{align}
M_0 -3t^{\frac 12+\epsilon}-(N^{\epsilon} \vee 2t^2/N)+c t N^{-\epsilon} \geq M_0 - N^{2\epsilon} + \theta t N^{-\epsilon}
\end{align} 
for, say, $\theta=c/2$; this implies by a union bound and~\eqref{eq:rw-domination}, that for every $t\leq T$,
\begin{equation*}
\mathbb P_{\cJ}(M_t\in \cB_{\epsilon,t,\theta}^c)\leq \bP(\cM_t - |\cR_t| \in \cA_{\epsilon,t,\theta}^c) \leq te^{-dt^{2\epsilon^2}}\,. \qedhere
\end{equation*}
\end{proof}

\subsection{Long time dynamics}\label{sec:long-time}
Using the bounds in \S\ref{sec:short-time} on the zero-temperature dynamics Markov chain, we can deduce the following conditions at time $T=N^{\frac 12+3\epsilon}$. 

\begin{proposition}\label{prop:everyone-positive-field}
If $M_0 = N^{\frac 12-\epsilon}$ and $T=N^{\frac 12+3\epsilon}$, there exists $\theta(p)>0$ such that 
\begin{align}\label{eq:everyone-positive-field}
\lim_{N\to\infty} \mathbb P \bigg(\bigcap_{i\notin \{Y_k\}_{k<T}} \big\{ m_i(T)>\epsilon p \theta N^{\frac 12 +2\epsilon}\big\}\bigg) =1\,.
\end{align}
\end{proposition}
\begin{proof}
We consider a fixed $i\notin \{Y_k\}_{k<T}$ and
prove the proposition using the decomposition of $m_{i}$ used in the proof of Proposition~\ref{prop:probability-positive-field} and then union bound over all such $i\notin \{Y_k\}_{k<T}$. First note that by a union bound with Lemmas~\ref{lem:good-update-sequence}--\ref{lem:coupling-concentration} and Proposition~\ref{prop:magnetization-bound}, if $\theta$ is as in Proposition~\ref{prop:magnetization-bound}, there exists $c(p)>0$ so that
\begin{align}
\mathbb P (M_T\in \mathcal B_{\epsilon,T-1,\theta},\Gamma_{\omega},\Gamma_{\cJ,T})\geq 1-O(Te^{-cN^{2\epsilon^2}})\,.
\end{align}
By a union bound, it suffices to prove that there exists $c>0$ such that for every $i\notin \{Y_k\}_{k<T}$, there exists $\theta>0$ such that
\begin{align}
\mathbb P_{\cJ} \left(m_i(T)\leq \epsilon p\theta N^{\frac 12+2\epsilon}  \mid \mathcal F_{T-1},\Gamma_{\omega},\Gamma_{\cJ,T}, \cB_{\epsilon,T-1,\theta}\right) \leq e^{-cN^{\epsilon}}\,.
\end{align}
Now suppose we are on the intersection of $\mathcal B_{\epsilon,T-1,\theta}, \Gamma_{\omega}$ and $\Gamma_{\cJ,T}$.
The couplings between $i$ and sites not in $\{Y_k\}_{k<T}$ are independent of $\mathcal F_{T-1}$ and are handled identically to the proof of Proposition~\ref{prop:probability-positive-field}, whence by Hoeffding's inequality applied to the difference of the binomial random variables in~\eqref{eq:binomials}, there exists $c(p,\theta)>0$ such that
\begin{align}\label{eq:contribution-havent-rung-yet}
 \mathbb P_{\cJ}\bigg(\sum_{j\notin \{Y_k\}_{k<T}} \sigma_j(0)J_{ij} \leq -\tfrac{p}{4}\theta N^{\frac 12+2\epsilon}\bigg) \leq e^{-cN^{\epsilon}}\,.
\end{align}
We now bound the contribution from couplings to sites $j\in \{Y_k\}_{k<T}-\cR_T$.
Under the event $\mathcal B_{\epsilon,T-1,\theta}$, the magnetization has $M_{T-1}\geq  M_0+ \theta N^{\frac 12 +2\epsilon}-N^{2\epsilon}$. Recalling the sets $Z_{++}(T),Z_{+-}(T),Z_{-+}(T),Z_{--}(T)$, as before, because $\Gamma_{\omega}$ holds, we have
\begin{align}\label{eq:difference-z-T}
\big||Z_{++}(T)| + |Z_{+-}(T)|-|Z_{--}(T)|-|Z_{-+}(T)|\big| & \leq N^{\frac 14 +3\epsilon}+N^{7\epsilon}\qquad\mbox{and} \nonumber \\
\qquad |Z_{-+}(T)|-|Z_{+-}(T)| & \geq \theta N^{\frac 12+2\epsilon} -N^{2\epsilon}-N^{7\epsilon}\,,
\end{align}
for some $\theta(p)>0$, implying that
\begin{align}
|Z_{++}(T)|+|Z_{-+}(T)|-|Z_{--}(T)|-|Z_{+-}(T)| \geq \theta N^{\frac 12+2\epsilon}-o(N^{\frac 12})\,.
\end{align}
In order to reveal jointly the couplings $\{J_{Y_k i}\}_{k<T}$ conditional on the history of the chain $(Y_k,S_{Y_k}(k))_{k<T}$, we use the revealing procedure defined in Definition~\ref{def:revealing-scheme}, taking $i$ to be $Y_T$; we know by Proposition~\ref{prop:condition-all-satisfactions} that under this revealing procedure, for every realization of $(S_{Y_k})_{k<T}$, the joint distribution of $\{J_{l i}\}_{l\in Z_{++}(T)}$ is dominated by a product measure of $\mbox{Ber}(p+N^{-\frac 12+4\epsilon})$ and dominates independent $\mbox{Ber}(p-N^{-\frac 12+4\epsilon})$ and the same holds for $\{J_{l i}\}_{l\in Z_{+-}(T)}$, and likewise when $l\in Z_{-+}(T)$ and $l\in Z_{--}(T)$. 
Then letting $Z_+=Z_{++}(T)+Z_{-+}(T)$ and $Z_-=Z_{--}(T)+Z_{+-}(T)$, by Hoeffding's inequality,
\begin{align}
\mathbb P_{\cJ}\bigg(\sum_{j\in \{Y_k\}_{k<T}-\cR_T} J_{ij} \sigma_j(T)  \leq \tfrac {p}2 \theta & N^{\frac 12+2\epsilon}\given \cF_{T-1}\bigg) \nonumber \\
& \leq \mathbb P_{\cJ} \bigg(\sum_{{j\in Z_+}} J_{ij} \leq p|Z_+|- \tfrac p4 \theta N^{\frac 12+2\epsilon}\given \cF_{T-1}\bigg) \nonumber \\
& \,\,\,\,\,\,+ \mathbb P_{\cJ} \bigg(\sum_{{j\in Z_-}} J_{ij} \geq p|Z_-| + \tfrac p4 \theta N^{\frac 12 +2\epsilon} \given \cF_{T-1}\bigg) \nonumber \\
 & \leq \,\, 2e^{-c N^{\epsilon}}\,,
\end{align}
for some $c(p,\theta)>0$.
Together with the bound of $N^{7\epsilon}$ on $\sum_{j\in \mathcal R_T} J_{ij} \sigma_j(T)$,  and a union bound, we obtain for some $c(p,\theta)>0$, under $\cB_{\epsilon,T-1,\theta}\cap \Gamma_\omega\cap \Gamma_{\cJ,T}$,
\begin{equation*}
\mathbb P_{\cJ} \bigg(\bigcup_{i\notin \{Y_k\}_{k<T}} \big\{m_i(T)\leq \tfrac p6 \theta N^{\frac 12+2\epsilon}\big\} \bigg) \lesssim   N e^{-cN^{\epsilon}}\,. \qedhere
\end{equation*}
\end{proof}

\begin{proposition}\label{prop:already-rung-pos-field}
If $M_0=N^{\frac 12-\epsilon}$ and $T'= N^{\frac 23}$, we have 
\begin{align}
\lim_{N\to\infty} \mathbb P \bigg(\bigcap_{i=1}^N \big\{ m_i(T')>0\big\}\bigg) =1\,.
\end{align}
\end{proposition}
\begin{proof}
First, consider the update sequence $Y_{T+1},...,Y_{T'}$. Dominating the number of updates there that are in $\{Y_k\}_{k\leq T}$ by $\mbox{Bin}(T',T/N)$, we see that the probability of that being at most $N^{\frac 12}$ is $1-O(e^{-c\sqrt n})$. At the same time, since $T'=o(N)$ and $M_0=N^{\frac 12-\epsilon}$, with probability at least $1-O(e^{-cN^{2/3}})$ there are at least $N^{\frac 23-\epsilon}$ distinct sites $i\in \{Y_{k}\}_{k=T+1}^{T'}-\{Y_k\}_{k\leq T}$ that have $\sigma_i(0)=-1$. Since both of these happen with $\mathbb P_{\omega}$-probability going to $1$ as $N\to\infty$ and are independent of $\cJ$ and $\cF_T$, suppose we are on the intersection of these events (and also $\Gamma_{\omega}$) and fix any such an update sequence. Then let $\cY_1,...,\cY_{n}\in \{Y_{k}\}_{k=T+1}^{T'}-\mathcal \{Y_k\}_{k\leq T}$ be the sequence of updates in $Y_{T+1},...,Y_{T'}$ with initial spin value $-1$ and not in $\cR_T$ (observe that $n\geq N^{\frac 23-\epsilon}$).

Assume also, by Proposition~\ref{prop:everyone-positive-field} that every site $i\notin \{Y_k\}_{k<T}$ has $m_i(T)\geq \epsilon p \theta N^{\frac 12 +\epsilon}$, as this occurs w.h.p. Now consider sites $i\in \{Y_k\}_{k<T}$. We will need the following lower bound on $m_i(T)$: if $T=N^{\frac 12 +3\epsilon}$ and $\Gamma_\omega$ holds,  there exists $c(p,\theta)>0$ such that
\begin{align}\label{eq:already-rung-field-lower-bound}
\lim_{N\to\infty} \mathbb P_\cJ \bigg (  \bigcup_{i \in \{Y_k\}_{k\leq T}} \big\{m_i(T) \leq -N^{\frac 12 + 4\epsilon} \big\} \bigg) \lesssim Te^{-cN^{2\epsilon}}\,.
\end{align}
We now define a new $\mathbb P_{\cJ}$ event $\Gamma^2_{\cJ}$ as 
\begin{align}
\Gamma^2_{\cJ} =\bigcap_{i\in \{Y_k\}_{k\leq T}}\bigg\{\sum_{j\in \cY_1,...,\cY_m} J_{ij}\geq (p-\epsilon)m,\, \Big|\sum_{j\notin \{Y_k\}_{k\leq T}} \sigma_{j}(0)J_{ij}\Big| \leq N^{\frac 12+\epsilon} \bigg\}\,.
\end{align} 
By a union bound and standard applications of Hoeffding's inequality as done before, along with the fact that $\Gamma^2_{\omega,T}$ holds, we see that for some $c>0$,
\begin{align}
\mathbb P_{\cJ} ((\Gamma_{\cJ}^2)^c) \leq 2Te^{-cN^{2\epsilon}}\,,
\end{align}
so without any loss, we also assume we are on the event $\Gamma_{\cJ}^2$. 

Recall that for every $j\notin \{Y_k\}_{k\leq T}$, $\sigma_j(0)=\sigma_j(T)$.
At the same time, we can assume a worst case bound on the couplings between site $i$ and other sites $j\in\{Y_k\}_{k\leq T}$, which is to say we take every such $J_{ij}\sigma_j(T)=-1$, contributing at most $-N^{\frac 12 +3\epsilon}$ to $m_i(T)$. Putting these together, we see that for update sequences in $\Gamma_\omega$, under the event $\Gamma_{\cJ}^2$ we have deterministically that for every $i\in \{Y_k\}_{k\leq T}$,
\begin{align}
m_i(T) \geq -N^{\frac 12+\epsilon}-N^{\frac 12+3\epsilon} =-o(N^{\frac 12+4\epsilon})
\end{align}
so that~\eqref{eq:already-rung-field-lower-bound} holds. Under the intersection of the events in~\eqref{eq:everyone-positive-field} and~\eqref{eq:already-rung-field-lower-bound}, we claim that by time $T'$, deterministically, every site will have positive effective field. Note that our update sequence $Y_{T+1},...,Y_{T'}$ is such that the dynamics only selects at most $N^{\frac 12}$ sites in $\{Y_k\}_{k\leq T}$ between times $T+1$ and $T'$ and as a result, by~\eqref{eq:everyone-positive-field}, for every $i\notin \{Y_k\}_{k\leq T}$, for every $T+1\leq t \leq T'$, its field satisfies 
\begin{align}
m_i(t)\geq m_i(T) -N^{\frac 12}\qquad \mbox{and therefore} \qquad m_i(t)>0
\end{align}
(only sites with nonpositive field can flip from $+1$ to $-1$ and decrease the field on $i$). At the same time, every time an update on a site $\cY_1,...,\cY_n$ occurs, that spin has positive field by the above, and so it flips from $-1$ to $+1$. Then since $n\geq N^{\frac 23-\epsilon}$ and $\Gamma_{\cJ}^2$ holds, for every $i\in \{Y_k\}_{k\leq T}$, its field satisfies
\begin{align}
m_i(T')\geq m_i(T)-N^{\frac 12}+ (p-\epsilon)N^{\frac 23-\epsilon} \qquad \mbox{and therefore} \qquad m_i(T')>0\,,
\end{align}
whenever $\epsilon$ is sufficiently small, concluding the proof of the proposition. 
\end{proof}

\begin{proof}[\textbf{\emph{Proof of Theorem~\ref{mainthm:1}}}]
By Proposition~\ref{prop:already-rung-pos-field}, by time $T'=N^{\frac 23}$, with high probability, every site $i$ has positive field. In that case by attractivity of the dynamics, whenever a negative site is selected to be updated, it flips to plus and those are the only possible spin flips, so that by the time every site has been updated again after time $T'$, the zero-temperature dynamics will have absorbed into the all-plus state.

Putting this together with the grand coupling of chains starting from every possible $M_0$ implies that if $\epsilon>0$ is sufficiently small, for every configuration with $M_0\geq N^{\frac 12 -\epsilon}$, 
\[\lim_{N\to\infty} \mathbb P ( M_\infty=N)=\lim_{N\to\infty} \mathbb P (\lim_{t\to\infty} M_t =N)=1\,. \qedhere
\]
\end{proof}
\begin{proof}[\textbf{\emph{Proof of Corollary~\ref{cor:qd}}}]
By spin flip symmetry, when $M_0\leq -N^{\frac 12-\epsilon}$, with high $\mathbb P_{\cJ,\omega}$-probability, $M_\infty=-N$.
Thus, for every $i=1,...,N$,
\begin{align}
\mathbb E_{\sigma(0)} \left[(\mathbb E_{\omega,\mathcal J}[\sigma_1(\infty)])^2\right] & = \mathbb E_{\sigma(0)} \left[(\mathbb E_{\omega,\mathcal J}[\sigma_1(\infty)])^2(\boldsymbol 1\{|M_0|\geq N^{\frac 12-\epsilon}\}+\boldsymbol 1\{|M_0|< N^{\frac 12-\epsilon}\})\right]\, \nonumber \\
& \geq \mathbb P_{\sigma(0)}(|M_0| \geq N^{\frac 12-\epsilon})   \min_{\sigma(0):|M_0|\geq N^{\frac 12-\epsilon}} (\mathbb E_{\omega,\mathcal J}[\sigma_1(\infty)])^2 \,.
\end{align}
Upper bounding the left-hand side by $1$ and taking limits as $N\to\infty$ on both sides, we obtain by Fact~\ref{fact:init-config} and Theorem~\ref{mainthm:1} that
\begin{align}
\lim_{N\to\infty} \mathbb E_{\sigma(0)} \left[ (\mathbb E_{\omega,\mathcal J} [\sigma_1(\infty)])^2\right] = 1\,. 
\end{align}
By Jensen's inequality,   
\begin{align}
q_D(N)= \mathbb E_{\sigma(0)} \left [  \mathbb E_{\mathcal J} [(\mathbb E_{\omega}[\sigma_i(\infty)])^2] \right ] \geq \mathbb E_{\sigma(0)}[(\mathbb E_{\omega,\mathcal J} [ \sigma_1 (\infty)])^2]
\end{align}
and trivially upper bounding $q_D(N)\leq 1$ implies $\lim_{N\to\infty} q_D(N)=1$.
\end{proof}

\section{Different behavior with heavy tails}\label{sec:heavy-tail}

\subsection{Heavy-tailed disordered CW model}\label{sec:heavy-tailed}
In this section, we demonstrate that the behavior above is not a general consequence of a mean-field disordered ferromagnetic system. We show that when we consider heavy-tailed coupling distributions, a completely different picture, regarding the structure of local-minima and the probability of ending up in them, emerges. One intuition for this comes from the fact that if we had considered an Erd\H{o}s--R\'enyi graph with $p=\lambda/N$ for $\lambda>0$ fixed, the results of \S\ref{sec:heavy-tail} on the dilute CW model would no longer hold as the underlying random graph would have many disconnected clusters of influence, and the dynamics would absorb in one of the many possible independent assignments of $\{\pm 1\}$ to the disconnected clusters. The arguments here are straightforward but we provide them to emphasize the contrast to the situation with light-tailed couplings.

If $J_{ij}\sim \mu$ is non-negative, we say it is \emph{heavy-tailed} if 
\begin{align}
\mathbb P(J_{ij}\geq x)=x^{-\alpha}L(x)
\end{align} for some $0<\alpha<1$ and some slowly varying $L(x)$ (i.e., for every $a>0$, $L$ satisfies $\lim_{x\to\infty} L(ax)/L(x)=1$). 
As before, for $\sigma\in\{\pm1\}^N$ define the Hamiltonian of this \emph{disordered Curie--Weiss model} by $H(\sigma)= -\frac 1N \sum_{i,j=1}^N J_{ij} \sigma_{i} \sigma_{j}$.
\begin{maintheorem}\label{thm:local-minima}
If $\{J_{ij}\}_{i,j}$ are non-negative heavy-tailed i.i.d.\ couplings with $0<\alpha<1$, then $H$ has at least one non-trivial local minimum on the hypercube $\{\pm 1\}^N$ with high probability. Moreover, the zero-temperature dynamics of the corresponding disordered Curie--Weiss model satisfies $\epsilon<q_D(N)<1-\epsilon$ for some $\epsilon>0$ uniform in $N$. 
\end{maintheorem}
\begin{remark}
One can also glean from the proof (specifically~\eqref{eq:expected-bully-bonds}) that on average, the number of local minima grows exponentially in $N$ as seen in many spin-glass models. While we expect that in the light-tailed setup, also with high probability there exist non-trivial local minima, there we guess the number in fact grows sub-exponentially.
\end{remark}
First recall the following classical theorems concerning sums and maxima of i.i.d.\ random variables with power-law tails (see e.g.,~\cite{Whittbook,Resbook}):

\begin{theorem}[Stable central limit theorem]\label{thm:stable-clt}
Let $Y_i$ be i.i.d.\ random variables satisfying the following conditions: there exists $0<\alpha<1$ such that $x^\alpha \mathbb P(|Y_1| \geq x)=L(x)$ is slowly varying and for some $\beta\in [-1,1]$, as $x\to\infty$,
\begin{align}
\mathbb P (Y_1 \geq x)/\mathbb P(|Y_1| \geq x) \longrightarrow (1+\beta)/2.
\end{align}
Then there exists a sequence $a_n$ given by $a_n^{-\alpha}L(a_n)=n^{-1}(\int_0^\infty x^{-\alpha} \sin x dx)^{-1}$ so that
\begin{align}
\frac {\sum_{i=1}^n Y_i}{a_n} \longrightarrow Z_{\alpha,\beta}\qquad \mbox{as $n\to\infty$}\,,
\end{align}
where $Z_{\alpha,\beta}$ is a standard $\alpha$-stable random variable with asymmetry parameter $\beta$.
\end{theorem} 

\begin{theorem}[Distribution of the maximum]\label{thm:max}
If $Y_1,...,Y_n$ are i.i.d.\ with $\mathbb P(Y_i\geq x)=x^{-\alpha} L(x)$ for $0<\alpha<1$ and $L(x)$ is slowly varying, then for each $x$, as $n\to\infty$,
\begin{align}
\mathbb P \big(\max_{i=1,..,n} Y_i  \leq b_n x \big) \longrightarrow \Phi_{\alpha}(x)=e^{-x^{-\alpha}}\boldsymbol1\{x> 0\}
\end{align}
where $b_n$ is the smallest sequence such that $b_n^{-\alpha} L(b_n) = n^{-1}$.
\end{theorem}

We prove that with high probability, the Hamiltonian $H$ has non-trivial local minima. In order to proceed we need to define bully bonds. 
For a coupling realization $\{J_{ij}\}_{1\leq i<j\leq N}$, the coupling $J_{ij}$ is a \emph{bully bond} if 
\begin{align}\label{eq:bully-bond}
J_{ij} > \max\bigg\{\sum _{k\neq j} |J_{ik}|,\sum_{k\neq i} |J_{jk}|\bigg\}\,,
\end{align}
and we define the event, $E_{ij}=\left\{J_{ij} >\max\{ \sum_{k\neq j} |J_{ik}|, \sum_{k\neq i} |J_{jk}|\}\right\}$.

\begin{proof}[\textbf{\emph{Proof of Theorem~\ref{thm:local-minima}}}]
We begin by proving that with high probability, there exist non-trivial local minima. Notice that if, with probability going to $1$ as $N\to\infty$, there exist at least two bully bonds, then $H$ has non-trivial local minima: if $E_{ij}, E_{kl}$ hold for $k,l\notin \{i,j\}$, both $\sigma_i=\sigma_j=+1, \sigma_k=\sigma_l=-1$ and $\sigma_i=\sigma_j=-1,\sigma_k=\sigma_l=+1$ are satisfied independent of $\sigma \restriction_{V-\{i,j,k,l\}}$ and neither combination is possible in a ground state of $H$.
By symmetry between $1\leq i<j \leq N/2$ and $N/2<k<l\leq N$, and a union bound, it suffices to prove that 
\begin{align}
\lim_{N\to\infty} \mathbb P_{\cJ}\bigg(\sum_{i,j=1}^{N/2} \boldsymbol 1\{E_{ij}\}>0\bigg)=1\,.
\end{align}
We show this by considering the probability that $\max_{1\leq i<j\leq N/2} J_{ij}$ is a bully bond: 
\begin{align}
\mathbb P_{\cJ}\bigg(J_{12} >  \max\{  \sum_{k\neq 1} J_{1k}, & \sum_{k\neq 2} J_{2k}\} \given J_{12}=\max_{1\leq i<j\leq N/2} J_{ij}\bigg) \geq \nonumber \\ 
& 1-\mathbb P_{\cJ}\bigg(\max_{1\leq i<j\leq N/2} J_{ij} \leq N^{\frac 2\alpha-\epsilon}\bigg) - 2\mathbb P_{\cJ}\bigg(\sum_{k=3}^{N} J_{1k}  \geq N^{\frac 2\alpha-\epsilon}\bigg)
\end{align}
for some small $\epsilon(\alpha)>0$. In the inequality, we used a union bound and the fact that the coupling distribution of $J_{1k}$ is decreased by conditioning on $J_{1k}$ not being the maximum coupling in a set.
Observe that the first probability above is $o(1)$ by Theorem~\ref{thm:max} and the fact that it is a maximum over order $N^2$ i.i.d.\ heavy-tailed random variables, and the second probability is $o(1)$ by Theorem~\ref{thm:stable-clt} and the fact that the sum is over $O(N)$ i.i.d.\ heavy-tailed random variables. By symmetry this implies that with high probability the maximum over all couplings $J_{ij}$ for $1\leq i<j \leq N/2$ is a bully bond, implying that with high probability there exist non-trivial local minima.  

It remains to prove that $q_D(N)$ is uniformly bounded away from $0$ and $1$. To do so, we begin by computing the expected number of bully bonds, 
\begin{align}\label{eq:expected-bully-bonds}
\mathbb E_{\cJ}\bigg[\sum_{ij} \boldsymbol 1\{E_{ij}\}\bigg]= \frac {N(N-1)}2 \mathbb P_{\cJ}(E_{12})\,.
\end{align}
Let $Z_{\alpha}=Z_{\alpha,\beta}$ with $\beta=1$ and let the sequence $a_N$ be defined as in Theorem~\ref{thm:stable-clt}; by definition of $a_N$ and independence,
\begin{align}
\mathbb P_{\cJ}(E_{12})  & \,\geq  \mathbb P_{\cJ} ( J_{12} > a_{2N})\mathbb P_{\cJ}\bigg(\sum_{k\neq 1} J_{1k}+\sum_{k\neq 2} J_{2k} \leq  a_{2N-2}\bigg) \nonumber \\
& \, \geq \frac{c}{2N}\left(\mathbb P_{\cJ}(Z_\alpha \leq 1)+o(1)\right)\geq \frac {c'}N\,,
\end{align}
for constants $c(\alpha),c'(\alpha)>0$. Note that this implies that $\mathbb E_{\cJ}[\sum_{ij} \boldsymbol 1\{E_{ij}\}] \geq \rho N$ for some $\rho(\alpha)>0$. Moreover, $\sum_{ij} \boldsymbol 1\{E_{ij}\}$ is bounded above by $N/2$ since every vertex can be adjacent to at most one bully bond. As a consequence, we have that 
\begin{align}
\mathbb P_{\cJ} \bigg(\sum_{1\leq i<j\leq N} \boldsymbol 1\{E_{ij}\} \geq \rho N\bigg) \geq 2 \rho\,.
\end{align}
We claim that this implies $q_D(N)$ is uniformly bounded away from $0$ and $1$. If with probability greater than $\epsilon>0$, there are at least $\delta N$ sites adjacent to bully bonds, then the probability of a single site $i$ being adjacent a bully bond $J_{ij}$ is at least $\epsilon \delta$. The contribution to $q_D$ on that event is $\frac 12$; this is because with $\mathbb P_{\sigma(0),\mathcal J}$-probability $\frac 12$, that bully bond is initially satisfied ($\sigma_i(0)=\sigma_j(0)$) and otherwise, it is equally likely that the dynamics absorbs with $\sigma_1=\sigma_j=+1$ as $\sigma_1=\sigma_j=-1$---this is completely determined by whether $\sigma_i$ or $\sigma_j$ is updated first. The definition of $q_D(N)$ then implies that $\frac{\epsilon \delta}2<q_D(N)<1-\frac{\epsilon \delta}2$ as desired.
\end{proof}

\subsection{Heavy--tailed spin glasses}
Using the above proof, one can derive similar estimates for heavy--tailed spin glasses, where $J_{ij}$ is now distributed as a symmetric heavy-tailed random variable with $\mathbb P(|J_{ij}|\geq x)=x^{-\alpha} L(x)$ for $L$ slowly varying and $0<\alpha<1$. The proof carries through as in the proof of Theorem~\ref{thm:local-minima}, with the sums becoming sums of absolute values of such random variables, which are in the same basin of attraction of fully asymmetric $\alpha$-stable random variables. 
Therefore, a straightforward adaptation of the proof of Theorem~\ref{thm:local-minima} shows that for such coupling distributions, the heavy-tailed spin glass has local minima with high probability, and its zero-temperature dynamical order parameter is uniformly bounded away from $0$ and $1$.

\bibliographystyle{abbrv}
\bibliography{Dilute-CW-New-proof}

\end{document}